\documentclass[12pt]{amsart}
\pdfoutput=1
\usepackage{amssymb,times,a4wide,esint,graphicx}

\theoremstyle{plain}
\newtheorem{theorem}{Theorem}[section]
\newtheorem*{theorem*}{Theorem}
\newtheorem{proposition}[theorem]{Proposition}
\newtheorem{corollary}[theorem]{Corollary}
\newtheorem{lemma}[theorem]{Lemma}

\theoremstyle{definition}
\newtheorem{definition}[theorem]{Definition}

\theoremstyle{remark}
\newtheorem*{remark}{Remark}

\newcommand{\supp}{{\rm supp}}
\newcommand{\R}{\mathbb{R}}
\renewcommand{\S}{\mathbb{S}}

\begin{document}
\title{Equivalent norms for polynomials on the sphere}

\author {Jordi Marzo}
\address{Departament de Matem\`atica Aplicada i An\`alisi
\newline \indent
Universitat de Barcelona, Gran via 585, 08071 Barcelona, Spain}
\email{jmarzo@mat.ub.es}

\author{Joaquim Ortega-Cerd\`a}
\address{Departament de Matem\` atica Applicada i An\`alisi
\newline \indent
Universitat de Barcelona, Gran Via 585, 08007-Barcelona, Spain}
\email{jortega@ub.edu}

\thanks{Supported by projects MTM2005-08984-C02-02 and 2005SGR00611}
\date{\today}
\begin{abstract}
    We find necessary and sufficient conditions for a sequence of sets
    $E_{L}\subset \S^{d}$
    in order to obtain the inequality
    \[\int_{\S^{d}}|Q_{L}|^{p} d\mu \le C_p \int_{E_{L}}|Q_{L}|^{p} d\mu,\quad
\forall L\ge 0,\]
    where $1\le p<+\infty$, $Q_L$ is any polynomial of degree smaller or equal
    than $L$, $\mu$ is a doubling measure and
    the constant $C_p$  is independent of $L$.
    From this description it follows an uncertainty principle for
    functions in $L^2(\S^{d})$.
    We consider also weighted uniform versions of this result.
\end{abstract}

 \maketitle

\section{Introduction}\label{SecPrelim}

The classical Logvinenko-Sereda theorem describes some equivalent norms
for functions in the Paley-Wiener space $PW^p_\Omega$, i.e. functions in
$L^p(\R^d)$ whose Fourier transform is supported in a prefixed bounded set
$\Omega\subset \R^d$.
\begin{theorem*}[Logvinenko-Sereda]
Let $\Omega$ be a bounded set and let $1\le p<+\infty$. A set $E\subset
\R^d$ satisfies
\[
 \int_{\R^d} |f(x)|^p\, dx \le C_p \int_{E} |f(x)|^p\, dx,\quad \forall
f\in PW^p_\Omega,
\]
if and only if
there is a cube $K\subset \R^{d}$  such that
    \[\inf_{x\in \R^{d}} |(K+x)\cap E|>0 .\]
\end{theorem*}
For a proof see \cite[pp.~112--116]{HJ} or the original \cite{LS}.

    Comparison norms results of this kind are known in other contexts, see
    \cite{HJ} and references therein for further information. The purpose of the
present paper is to prove similar comparison results for $L^{p}$
norms of polynomials on the unit sphere $\S^{d}$.

In what follows $\sigma$ will denote the surface measure in $\S^{d}$.
We will prove the following theorem:

\begin{theorem*}
    Let $1\le p <\infty$.    A sequence of sets $\mathcal{E}=\{ E_{L}
\}_{L\ge 0}$ in $\S^{d}$ satisfies
    \begin{equation}\label{L2}
    \int_{\S^{d}}|Q_{L}|^{p} d\sigma \le C_p
    \int_{E_{L}}|Q_{L}|^{p} d\sigma,\quad \forall L\ge 0,
    \end{equation}
    where $Q_L$ is any polynomial of degree smaller or equal than  $L$ and
    the constant $C_p$  is independent of $L$ if and only if
    \[ \inf_{L\in\mathbb{N},1-|z|=1/L} h_z(E_L)>0,\]
    where $h_z(F)$ is the harmonic extension of $\chi_F$ to a point
$z$ in the interior of the ball.
\end{theorem*}
    A more general (and precise) version will be stated and proved later on, see
Theorem~\ref{main-theorem} once we have introduced some definitions and
notation.

    From this Theorem it follows an uncertainty
    principle for functions in $L^{2}(\S^{d})$.
    For any $f\in L^{2}(\S^{d})$ we have the spherical harmonics expansion
    $f=\sum_{\ell \ge 0}P_\ell(f)$,
    where $P_\ell$ is the orthogonal projection from $L^2(\S^d)$ to
$\mathcal{H}_\ell$

\begin{corollary}
    For a  set $E\subset \S^{d}$ let
    $\delta=\inf_{1-|z|=1/L}h_z(E)$.
    There exists a constant $C>0$ depending only on $\delta$
    such that for any $f\in L^{2}(\S^{d})$
\begin{equation}
       \label{strong}
    \int_{\S^{d}}|f(u)|^{2} d\sigma(u)\le C\left(
    \int_{E}|f(u)|^{2} d\sigma(u)+\sum_{\ell
>L}\|P_\ell (f)\|^{2} \right).
\end{equation}
\end{corollary}
   The proof of the Corollary amounts to show that \eqref{strong} is
equivalent to the inequality \eqref{L2} and it can be found in
\cite[3.1.1.A)~pp.~88--89]{HJ}.

\subsection{Preliminaries and Statements}

    In $\S^{d}$ we take the geodesic distance
    \[d(u,v)=\arccos \langle u, v\rangle ,\quad u,v\in \S^{d},\]
    and let $B(\omega , \delta ) \subset \S^{d}$ denote
    the geodesic ball of center $\omega\in \S^{d}$ and
    radius $\delta>0$.

    Let ${\mathcal{H}}_{\ell}$ be the spherical harmonics of degree $\ell$
    i.e. the restrictions to the unit sphere
    $\S^{d}$ of the homogeneous harmonic polynomials in $d+1$ variables of
    degree $\ell$.
    Let
    $\Pi_{L}=\bigcup_{\ell=0}^{L}{\mathcal{H}}_{\ell}$ denote
    the spherical harmonics of degree less or equal than $L$.
    Observe that the restriction to $\S^{d}$ of any polynomial in $d+1$
    variables of degree $\le L$ belongs to
    $\Pi_{L}$.

    In the Hilbert space $L^{2}(\sigma)$
    let us denote by
    $Y_{\ell}^{1},\dots ,Y_{\ell}^{h_{\ell}}$ an orthonormal basis of
${\mathcal{H}}_{\ell}$.
    Taking all these bases for $\ell=0,\dots L$ together we get an orthonormal
basis for $\Pi_{L}$.

    It is well known that the reproducing kernel for $\Pi_{L}$ is

\[ K_{L}(u,v)=\sum_{\ell=0}^{L}\sum_{j=1}^{h_{\ell}}Y_{\ell}^{j}(u)\overline{Y_{
\ell}^{j}(v)},\quad u,v\in \S^{d},\]
    and this expression does not depend on the choice of the bases. Using the
    Christoffel-Darboux formula (see for instance \cite{J06}) we obtain
    \[ K_{L}(u,v)=\frac{\kappa_{d,L}}{\sigma(\S^{d})}P_{L}^{(d/2,d/2-1)}
    (\langle u, v\rangle),\]
    where
    $P_{L}^{(\alpha,\beta)}$ stands for the Jacobi polynomial of degree $L$ and
    index $(\alpha,\beta)$
    and
    $\kappa_{d,L}\footnote{Here and in
    what follows $\sim$ means that the ratio of the two sides is bounded from
above and from
    below by two positive constants.}\sim L^{d/2}$, as $L\to \infty$.
    From now on we denote $\lambda=(d-2)/2$.

    Finally, we recall an estimate, \cite[p.~198]{S}, that will be used later
on:
\begin{equation}
       \label{estimate}
    P_{L}^{(1+\lambda,\lambda)}(\cos
    \theta)=\frac{k(\theta)}{\sqrt{L}}\left\{ \cos
\left((L+\lambda+1)\theta-\frac{(d+1)\pi}{4} \right)+\frac{O(1)}{L\sin
    \theta}\right\},
\end{equation}
    if $c/L\le \theta \le \pi-(c/L)$, where
    \[k(\theta)=\pi^{-1/2} \left( \sin
\frac{\theta}{2}\right)^{-\lambda-3/2}\left( \cos
    \frac{\theta}{2}\right)^{-\lambda-1/2}.\]

\begin{definition}
    We say that a measure $\mu$
    is \emph{doubling} if there exist a constant $C>0$ such that for any $u\in
\S^{d}$ and any $\delta>0$,
    \[ \mu(B(u,2\delta))\le C\mu(B(u,\delta)),\]
    For such a measure $\sup_{u,\delta} \mu(B(u,2\delta))/\mu(B(u,\delta))$
    is called the doubling constant of $\mu$.
\end{definition}

    It can be seen
(see for instance \cite[Lemma 2.1.]{MT})
    that for $\mu$ doubling there exist a $\gamma>0$ such that for $r>r'$
\begin{equation}
                   \label{doubling}
    \left( \frac{r}{r'}\right)^{1/\gamma} \lesssim
\frac{\mu(B(u,r))}{\mu(B(u,r'))}
    \lesssim \left( \frac{r}{r'}\right)^{\gamma},
\end{equation}
    with constants depending only on the doubling constant of $\mu$.

    Mimicking the Euclidean situation we define the following concept.

\begin{definition}
    Let $1\le p< \infty$ and let $\mu$ be a doubling measure.
    We say that the sequence of sets $\mathcal{E}=\{ E_{L} \}_{L\ge 0}\subset
\S^{d}$ is
    $L^{p}(\mu)$-\emph{Logvinenko-Sereda}
    if there exists a constant $C_{p}>0$ such that for any $Q\in \Pi_{L}$ and
any $L$
\begin{equation}
           \label{defineq}
    \int_{\S^{d}}|Q(u)|^{p}d\mu (u) \le C_{p}\int_{E_{L}}|Q(u)|^{p}d\mu (u).
\end{equation}
\end{definition}

\begin{definition}
    The sequence of sets $\mathcal{E}=\{ E_{L} \}_{L\ge 0}\subset \S^{d}$ is
$\mu$-\emph{relatively dense}
    if
    there exist $r>0$ and $\varrho>0$ such that
\begin{equation}
      \label{rd}
    \inf_{u \in \S^{d}}\frac{\mu(E_{L}\cap B(u,r/L))}{\mu(B(u,r/L))}\ge
\varrho>0,
\end{equation}
    for all $L$. When $\mu$ is the Lebesgue measure we say
    that $\mathcal{E}$ is relatively dense.
\end{definition}

    Now we can state our main result.

\begin{theorem}
\label{main-theorem}
    Let $\mathcal{E}=\{ E_{L} \}_{L \ge 0}$ be a sequence of sets in $\S^{d}$.
    $\mathcal{E}$ is $L^{p}(\mu)$-Logvinenko-Sereda for some
    $1\le p< \infty$ and $\mu$ a doubling measure if and only if $\mathcal{E}$
    is $\mu$-relatively dense.
\end{theorem}

If $\mu$ is absolutely continous with an $A_\infty$ weight it is possible to
reformulate the $\mu$-relatively density in terms of the harmonic extension.

   For a weight $\omega\ge 0$ in $\S^{d}$
    we denote
    \[\omega(E)=\int_{E}\omega(u)   d\sigma(u),\quad E\subset \S^{d}.\]

\begin{definition}
   \label{ainfty}
   A weight $\omega$ belongs to $A_{\infty}$ if there exist constants
$B,\beta>0$
    such that
    for any $E\subset B(u,\delta)$
    measurable
\begin{equation}
           \label{defainfty}
    \omega(B(u,\delta))\le B
    \left(\frac{\sigma(B(u,\delta))}{\sigma(E)}\right)^{\beta}\omega(E).
\end{equation}
\end{definition}
    It is well know that that an $A_{\infty}$ weight defines a doubling measure
but the converse is not true, see \cite{FM}.

\begin{remark} Due to the reversibility of condition \eqref{defainfty}, see
\cite[chap.~V,~1.7]{Stein},
    being relatively dense is equivalent to the same condition for the measure
defined with $\omega\in A_{\infty}$.
\end{remark}

    We recall that for $x\in \R^{d+1}$ with $|x|<1$ the harmonic measure
    of subset $F\subset \S^{d}$ with respect to $x$ is
\[h_{x}(F)=\frac{1}{\sigma(\S^{d})}\int_{F}\frac{1-|x|^{2}}{|x-u|^{d+1}}
 d\sigma(u)=
    \frac{1}{\sigma(\S^{d})}\int_{F}P(x,u) d\sigma(u),\]
    and $P(x,u)$ is the Poisson kernel in $\S^{d}$.
    The next result is a version for $\S^{d}$ of the one proved in
    \cite[p.~114]{HJ}.

\begin{lemma}
      \label{lemma}
    The sequence $\{E_{L}\}_{L\ge 0}\subset \S^{d}$
    is relatively dense if and only if there exists
    $\alpha>0$ such that
\[h_{x}(E_{L})\ge \alpha,\mbox{ for all } x\in
\R^{d+1}\mbox{ with }|x|=1-1/L.\]
\end{lemma}

\begin{proof}
    Observe that both conditions are rotation invariant.
    For $u$ such that $d(u,N)<r/L$ we have $C L^{d}\le P(|x|N,u)\le 2L^{d}$,
    where $C>0$ is a constant depending on $r$ and $d$.
    For $\theta=d(u,N)>r/L$
    \[P(|x|N,u)\lesssim
    \frac{\frac{2}{L}-\frac{1}{L}}{\sin^{d+1}\frac{\theta}{2}}\lesssim
    \frac{L^{d}}{r^{d+1}}.\]
    These bounds are all we need to prove the result.
    In one direction
    \[h_{|x|N}(E_{L})\gtrsim L^{d}\sigma(E_{L}\cap B(N,r/L))\gtrsim \varrho>0.\]
    Conversely
\begin{align*}
    \sigma(\S^{d}) \alpha & \le \int_{E_{L}} P(|x|N,u)  d\sigma(u)
    \le
    \int_{E_{L}}(\chi_{B(N,r/L)}(t)+ \chi_{B(N,r/L)^{c}}(t))
    P(|x|N,u) d\sigma(u)
    \\
    &
    \le
    2L^{d}\sigma( E_{L}\cap B(N,r/L))
    +
    \sum_{\log_{2}\pi L/r \ge j\ge
0}\int_{\frac{2^{j}r}{L}<d(u,N)<\frac{2^{j+1}r}{L}}
    \chi_{E_{L}}(u)P(|x|N,u) d\sigma(u)
    \\
    &
    \le
    C_{r}\frac{\sigma(E_{L}\cap B(N,r/L))}{\sigma(B(N,r/L))}+
    \frac{C}{r} \sum_{j\ge 0} \frac{1}{2^{dj}}.
\end{align*}
    Taking $r>0$ big enough we get the result.
\end{proof}

\begin{remark} We have proved that there exist $r,\varrho$ such that
    $\sigma (E_{L}\cap B(u,r/L))\ge \varrho \sigma(B(u,r/L))$ for  $L$ big enough
if and only if
    there exists $\alpha$ such that
    $h_{(1-1/L)u}(E_{L})\ge \alpha$. This new formulation depends only on one
parameter.
\end{remark}

    In Theorem~\ref{main-theorem} when the dimension $d=1$ there are already
some results known. In this case it is possible to replace polynomials by
holomorphic
polynomials. If moreover $\mu$ is the Lebegue measure the space of holomorphic
polynomials can be seen as a model space, so
    Volberg result \cite{V}
    extending the original theorem of Logvinenko and Sereda
    to model spaces apply.
    Also when $d=1$ and the measure $\mu$ is an $A_{\infty}$ weight,
    the sufficiency of condition \eqref{rd}
    was proved in \cite[Theorem~5.4]{MT}.

    Condition \eqref{rd} is true for some $\omega\in A_{\infty}$ if and only if
    it is true for the Lebesgue measure $\sigma$. So we have comparison of norms
for
any $\omega\in A_{\infty}$
    if and only if we have \eqref{rd} for the Lebesgue measure $\sigma$.
    The discussion following \cite[Theorem~5.4]{MT} shows that this is not true
for arbitrary
    doubling measures.

    To the best of our knowledge, for dimensions greater
than one Theorem~\ref{main-theorem} is
new, even in the case of Lebesgue measure.

    The outline of this paper is as follows.
    In Section~\ref{SecU}
    we will prove Theorem~\ref{main-theorem}.

    In Section~\ref{Secinfty} we deal with the uniform norm case.
    In this setting we have an analogous result to Proposition~\ref{prop2},
namely Theorem~\ref{lsinfty}.
    To consider weighted versions of this result an obvious requirement is to
take
    weights bounded above. We take the reverse H\"older class $RH_{\infty}$ of
those
    weights satisfying reverse H\"older inequalities in a uniform way.
    This class that was also introduced in \cite{MT} for the one-dimensional
case,
    is shown to be optimal in a certain sense.

\section{Main Results}\label{SecU}

\begin{proposition}
                           \label{prop1}
    Let $1\le p< \infty$, $\mu$ be a doubling measure
    and let $\mathcal{E}=\{ E_{L} \}_{L \ge 0}$ be a sequence of sets in
$\S^{d}$.
    If $\mathcal{E}$ is $L^{p}(\mu)$-Logvinenko-Sereda, then it is
$\mu$-relatively dense.
\end{proposition}

\begin{proof}
    We focus on $d\ge 2$ but only minor changes will prove the one-dimensional
case.
    The strategy is to apply the $L^{p}(\mu)$-comparison of norms to a power of
the reproducing kernel
    and to use classical estimates on the Jacobi polynomials.

    Let $Q(v)=(P_{L}^{(1+\lambda,\lambda)}(\langle v ,N \rangle))^{\ell}\in
\Pi_{\ell
L}$ and let $0<r<<R$. We have by hypothesis,
\begin{align}\label{carleson}
    \int_{B(N,r/L)} & |Q(v)|^{p}d\mu (v) \le \int_{\S^{d}}|Q(v)|^{p}d\mu (v)\le
C
    \int_{E_{\ell L}}|Q(v)|^{p}d\mu (v)
    \\ \notag
    &
    \lesssim
    \int_{E_{\ell L}\cap B(N,R/L)}|Q(v)|^{p}d\mu (v)+\int_{\S^{d}\setminus
B(N,R/L)}|Q(v)|^{p}d\mu (v).
\end{align}
    Observe that $Q$ reach its maximum in $N$, \cite{S} so applying Bernstein's
    inequality to the polynomial restricted to a great circle
    we get for any $v$ such that $d(v,N)<r/L$
    \[|Q(v)-Q(N)|\le |Q(N)|\ell r.\]
    Therefore for $r$ small enough  we have
$|Q(v)|^{p}\sim
|Q(N)|^{p}$ if $d(v,N)<r/L$.
    We can bound the  integral in the left hand side of \eqref{carleson} as
    \[\int_{B(N,r/L)}|Q(v)|^{p}d\mu (v)\gtrsim
(P_{L}^{1+\lambda,\lambda}(1))^{p\ell}\mu(B(N,r/L))
    \sim L^{\frac{p\ell d}{2}}\mu(B(N,r/L)).\]
    Since $| P^{(1+\lambda,\lambda)}_{L}(\cos \theta)|\lesssim L^{\lambda}$ for
$\pi-\frac{R}{L}\le \theta \le \pi$,
\begin{align*}
    L^{\frac{p \ell d}{2}}  \mu(B(N,r/L)) &\lesssim
    L^{\frac{p \ell d}{2}} \mu(E_{\ell L}\cap B(N,R/L))+
    L^{\frac{p \ell (d-2)}{2}} \mu(B(S,R/L))
    \\
    &
    +\int_{R/L<d(v,N)<\pi-R/L}|Q(v)|^{p}d\mu (v).
\end{align*}

To control the last integral we may use Szeg\"o  estimate \eqref{estimate}
\begin{align*}
     & \int_{R/L<d(v,N)<\pi-R/L} |Q(v)|^{p}d\mu (v)
    \lesssim
    L^{-p\ell/2}\int_{R/L<d(v,N)<\pi/2} \left|    \sin^{d+1} \frac{d(v,N)}{2}
     \right|^{-\ell p/2}
    d\mu(v)
    \\
    &
    +
    L^{-p\ell/2}\int_{\pi/2 <d(v,N)<\pi-R/L} \left|
    \cos^{d-1} \frac{d(v,N)}{2}      \right|^{-\ell p/2}
    d\mu(v)=I+II.
\end{align*}

    For part $I$ we take $\ell$ big enough to get $C(\mu)<2^{pl(d+1)/4}$, where
$C(\mu)$ is the doubling constant of $\mu$.
    We split the sphere in diadic ``bands'' around the north pole and using the
doubling property for $\mu$ we get
\begin{align*}
    L^{p\ell /2}I&
    \lesssim \int_{R/L<d(v,N)} \frac{1}{d(v,N)^{(d+1)lp/2}} d\mu (v)
    \le
    \sum_{J\ge j\ge 0} \int_{2^j R/L<d(v,N)<2^{j+1}R/L} \frac{d\mu
(v)}{(2^jR/L)^{(d+1)lp/2}}
\\
&\le
    \sum_{J\ge j\ge 0}  \frac{\mu(B(N,2^{j+1}R/L))}{(2^jR/L)^{(d+1)lp/2}}
    \le
    \frac{\mu(B(N,R/L))}{(R/L)^{(d+1)lp/2}} \sum_{j\ge 0}  \left(
\frac{C(\mu)}{2^{\alpha+\lambda}}\right)^{j}
    \lesssim \frac{\mu(B(N,R/L))}{(R/L)^{(d+1)lp/2}},
\end{align*}
    where ${\mathbb N} \ni J\ge \log_{2} (\pi L/R)$.

    For part $II$ the same computation taking dyadic ``bands'' around the south
pole shows that
\[
    L^{p\ell/2} II \lesssim \frac{\mu(B(S,R/L))}{(R/L)^{(d-1)lp/2}}
\lesssim \mu(B(S,R/L)) L^{(d-1)lp/2}.
\]
We use now property \eqref{doubling} and the $\gamma$ given there for
$\mu$ to estimate $\mu(B(S,R/L))$. If we put all estimates together and
for $\ell$ big enough we get
\begin{align*}
    \mu & (B(N,r/L))
    \lesssim
    \mu(E_{\ell L}\cap B(N,R/L))+
    L^{-p\ell} \mu(B(S,R/L))
    +\frac{\mu(B(N,R/L))}{R^{(d+1)lp/2}}
    \\
    &
    \lesssim
    \mu(E_{\ell L}\cap B(N,R/L))+
    \frac{R^{\gamma}}{L^{\gamma+p\ell}}+
    \left( \frac{R}{r} \right)^{\gamma}
\frac{\mu(B(N,r/L))}{R^{(d+1)lp/2}}.
\end{align*}
    As $\mu(B(N,r/L))\ge (r/L)^{1/\gamma}$ the second
    term is $o(\mu(B(N,r/L)))$ when $L\to \infty$
    for $\ell$ big enough. For the third term we choose $\ell$ such that
    $(R/r)^{\gamma}\le R^{(d+1)lp/2}/2$. Thus picking $\ell$ big enough we have
proved that
\[
 \mu (B(N,r/L)) \lesssim \mu(E_{\ell L}\cap B(N,R/L)),\quad \text{if }L\ge L_0
\]
Of course, the constants do not depend of the center of the balls being the
north pole.
Moreover by the doubling property $\mu (B(N,R/L))\simeq \mu (B(N,r/L))$. By
choosing a bigger $R$ we get
\[
 \mu (B(z,R/L)) \lesssim \mu(E_{\ell L}\cap B(z,R/L)),\quad \forall z\in
\S^d,\ L\ge 0.
\]
Finally we have only controlled the density of the sequence of
sets $\{E_{lL}\}_{L\ge 0}$. But we could have used the same
argument to the sequence $\mathcal{E}'=\{E_{L+1}\}_{L\ge 0}$ from
the very beginning and we will obtain then a control of the
density the sets $\{E_{lL+1}\}_{L\ge 0}$. By repeating the
argument $l$ times we get the desired result.
\end{proof}

\begin{remark} The somehow simpler polynomials
    \[\left( \frac{1+\langle v,N\rangle}{2}\right)^{L\ell}\mbox{ or }
    \left( \frac{1-\langle v, N\rangle^{L+1}}{(L+1)(1-\langle v,
N\rangle)}\right)^{\ell}\]
    that peak at $N$ and have been considered in other contexts
    do not decrease fast enough near the pole north to be chosen as test
functions for
the comparison of norms as we did with the polynomial $Q$ above.
\end{remark}

\begin{proposition}
                           \label{prop2}
    If $\{ E_{L} \}_{L\ge 0}$ is $\mu$-relatively dense for some doubling
measure $\mu$
    then it is $L^{p}(\mu)$-Logvinenko-Sereda for any
    $1\le p< \infty$.
\end{proposition}

\begin{proof}
    We consider a regularized version of $\mu$
    \[\mu_{L}(u)=\frac{\mu (B(u,1/L))}{\sigma (B(u,1/L))},\quad L\ge 0.\]
    By Corollary 3.4. in \cite{Dai} we have
    \[\int_{\S^{d}}|Q_{L}(u)|^{p}d\mu (u)\sim \int_{\S^{d}}|Q_{L}(u)|^{p}
\mu_{L}(u)  d\sigma (u),
    \quad Q_{L}\in \Pi_{L}.\]
    The regularization of $\mu$ is pointwise equivalent to a polynomial. Indeed,
    there exists $R_{L}\in \Pi_{L}$ nonnegative such that
    for any $u\in \S^{d}$
    \[\mu_{L}(u)\sim R_{L}(u)^{p},\]
    with constant depending only on $d$, the doubling constant for $\mu$ and
$p$, see \cite[Lemma~4.6]{Dai}.
    Given $Q_{L}\in \Pi_{L}$ let $M_{2L}\in \Pi_{2L}$ such that
$M_{2L}=Q_{L}R_{L}$ in $\S^{d}$.
    Following an idea of D.~H.~Luecking \cite{Lue83} we consider, for $\epsilon>0$ and $r>0,$ the set
    of points $z\in \S^{d}$
such that $M_{2L}(z)$
    has the same size as its average, i.e.
    \[A=A_{\epsilon,r,M_{2L}}
    =\left\{ z \in \S^{d}:|M_{2L}(z)|^{p}\ge \epsilon
\fint_{\mathbb{B}(z,r/L)}|M_{2L}(u)|^{p}dm(u)  \right\}.\]
    Most of the norm of $M_{2L}$ is concentrated on $A,$
\begin{align*}
    \int_{\S^{d}\setminus A} & |M_{2L}(z)|^{p} d\sigma (z) \le \epsilon
    \int_{\S^{d}\setminus A}\left(\fint_{\mathbb{B}(z,r/L)}|M_{2L}(u)|^{p}dm(u)\right)
d\sigma(z)
    \\
    &
    \le
    \epsilon \int_{|1-|u||<r/L}
    |M_{2L}(u)|^{p}
    \left( \int_{\S^{d}\setminus A}
\frac{\chi_{\mathbb{B}(z,r/L)}(u)}{m(\mathbb{B}(z,r/L))} d\sigma
    (z)\right)dm(u)
    \\
    &
    \lesssim \epsilon L
    \int_{|1-|u||<r/L}
    |M_{2L}(u)|^{p}
    dm(u)\sim \epsilon \int_{\S^{d}}|M_{2L}(z)|^{p} d\sigma (z),
\end{align*}
    using \cite[Corollary~4.3]{J06} in the last estimate, the constants are independent of $L.$

    Thus it is enough to
    show that
\begin{equation*}
    \int_{A}|M_{2L}(u)|^{p} d\sigma(u) \lesssim
\int_{E_{L}}|Q_{L}(u)|^{p}d\mu(u).
\end{equation*}

    All we need to prove is the existence of a constant $C>0$ such that for all
$\omega \in A$
\begin{equation} \label{main-des}
    |Q_{L}(\omega)|^{p}\le \frac{C}{\mu(B(\omega,r/L))}\int_{B(\omega,r/L)\cap
    E_{L}}|Q_{L}(u)|^{p}d\mu(u).
\end{equation}
    Indeed, if this is the case then
\begin{align*}
    \int_{A}|M_{2L}(\omega)|^{p} d\sigma (\omega)
    &
    \le
    C\int_{E_{L}}|Q_{L}(u)|^{p}\int_{\S^{d}}
    \frac{\chi_{B(\omega,r/L)}(u)}{\mu (B(\omega,r/L))}\mu_{L}(\omega) d\sigma
(\omega)
    d\mu(u)
    \\
    &
    \lesssim
    \int_{E_{L}}|Q(u)|^{p}
    d\mu (u).
\end{align*}

    To prove \eqref{main-des} we argue by contradiction. If
    (\ref{main-des}) is false there are
    for any $n\in {\mathbb N}$ polynomials
    $Q_{n}\in \Pi_{L_{n}}$ and $\omega_{n}\in A$ such that
\begin{equation*}
    |Q_{n}(\omega_{n})|^{p}>
\frac{n}{\mu(B(\omega_{n},r/L_{n}))}\int_{B(\omega_{n},r/L_{n})
    \cap E_{L_{n}}}|Q_{n}(u)|^{p}d\mu(u).
\end{equation*}

    Since $\mu$ is doubling then $R_{L_{n}}(\omega_{n})\sim R_{L_{n}}(u)$ for any $u\in
B(\omega_{n},r/L_{n}).$
    Let $M_{n}\in \Pi_{2L_{n}}$ be such that
    $M_{n}=Q_{n}R_{L_{n}}$ in $\S^{d}$
\begin{equation}
               \label{contra}
    |M_{n}(\omega_{n})|^{p}\gtrsim
\frac{n}{\mu(B(\omega_{n},r/L_{n}))}\int_{B(\omega_{n},r/L_{n})
    \cap E_{L_{n}}}|M_{n}(u)|^{p}d\mu(u).
\end{equation}

    By means of a rotation, a dilation and a translation we send
    $\omega_{n}$ to the origin in $\R^{d+1}$, the ball
$\mathbb{B}(\omega_{n},r/L_{n})$
    to $\mathbb{B}(0,1)\subset \R^{d+1}$ and the set $E_{L_{n}}$
    to
    \[E_{n}\subset \partial \mathbb{B}(-(L_{n}/r) N,L_{n}/r)\cap
    \mathbb{B}(0,1).\]
    The composition of these applications with our harmonic polynomials
    $M_{n}$ are harmonic functions $f_{n}$ that, after normalization, we can
    assume
    that satisfy
    \[\int_{\mathbb{B}(0,1)}|f_{n}|^{p}dm=1.\]
    The subharmonicity of $|f_{n}|^{p}$ and the fact that $\omega_{n}\in A$
    tells us that
\begin{equation*}
    \epsilon \lesssim |f_{n}(0)|^{p}\lesssim 1,
\end{equation*}
    and this property together with \eqref{contra} yields
\begin{equation}                                                \label{transf}
    \frac{1}{n}\gtrsim \int_{\mathbb{B}(0,1)\cap
E_{n}}|f_{n}(u)|^{p}d\mu_{n}(u),
\end{equation}
    where $\mu_{n}$ is the push forward of the measure $\mu,$ supported in
    $\partial \mathbb{B}(-(L_{n}/r) N,L_{n}/r)\cap
    \mathbb{B}(0,1),$ and normalized in such a way that $\mu_{n}(\mathbb{B}(0,1))=1.$

    We have that
    $\{ f_{n} \}$ is a normal family in $\mathbb{B}(0,1)$ and therefore
    there exist a subsequence that converges locally uniformly on $\mathbb{B}$
to
    an harmonic function that we call $f$.

    We observe that the relative density hypothesis yields
    \[\inf_{n}\mu_{n}(E_{n}\cap \mathbb{B}(0,1))>0.\]

    Let $\tau$ be a weak-$*$ limit of a subsequence of
    $\tau_{n}=\mu_{n}\chi_{E_{n}},$ having $\supp \tau\subset \R^{d}\times \{ 0 \}$ and
    $\tau\not\equiv 0.$
    We will consider the measure $\tau_{n}$ that has support in
    $\partial \mathbb{B}(-(L_{n}/r) N,L_{n}/r)\cap \mathbb{B}(0,1)$
    as having support in $\R^{d}\times \{ 0 \}$. To do so
    we define the measure $\tilde{\tau}_{n}$ as the ``projection'' of the
measure $\tau_{n}$ to
    $\R^{d}\times \{ 0 \}$, i.e.
    $\tilde{\tau}_{n}(A)=\tau_{n}(A\times [-1,1])$, for $A\subset \R^{d}$.

    We observe that $f$ restricted to $\R^{d}\times \{ 0 \}$ is real analytic.
    Condition \eqref{transf} implies that $f=0$ $\tau$-a.e.
    and therefore $\supp \tau\subset \{ f_{|\R^{d}\times \{0\}}=0 \}$.

    We want to show that
    $\supp \tau\subset \R^{d}$ (identifying $\R^{d}\times \{0\}$ and $\R^{d}$)
    cannot lie on a real analytic $(d-1)$-dimensional submanifold $S\subset
\R^{d}$
    (the worst case).
    We argue by contradiction. Let $x\in \supp \tau \subset S$ and $\delta>0$ be
such that
    $\tau(B(x,\delta))=\epsilon>0$.

    We can consider for any $y\in B(x,\delta)\cap S$ the unitary vector $\nu_{y}$ in $\R^{d}$ normal
    to $S$ in the point $y$ (see figure \ref{dibuix}) and define the ``square''
    $B(x,\delta) \subset R_{x}$
    \[R_{x}=\{ y+\eta\nu_{y}:y\in B(x,\delta)\cap S,|\eta |<\delta \}.\]

\begin{figure}
    \begin{center}
    \includegraphics[scale=0.7]{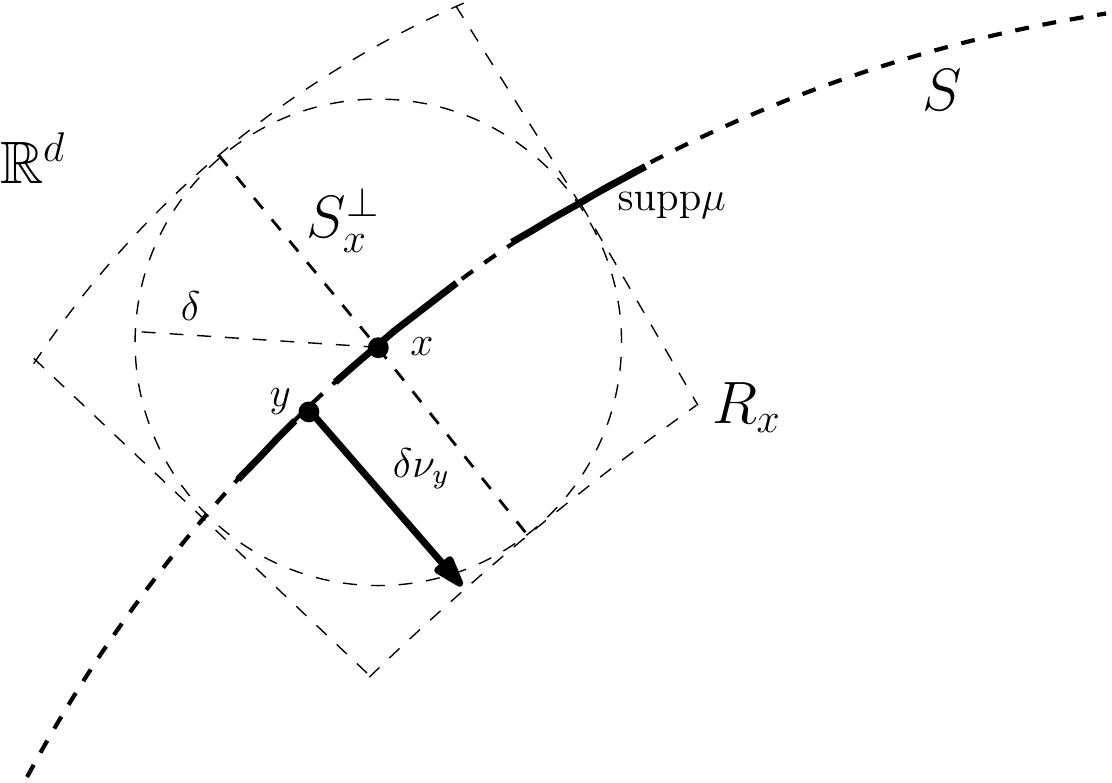}
    \end{center}
    \caption{ \label{dibuix}}
\end{figure}

    Now
    we can define measures $\nu_{n}$ in $S^{\perp}_{x}=\{ x+\eta\nu_{x}:|\eta
|<\delta\}$
    just by taking for $A\subset S_{x}^{\perp}$ the set
    $\tilde{A}\subset (-\delta,\delta)$ such that $x+\tilde{A} \nu_{x}=A$ and
defining
    \[\nu_{n}(A)=\tilde{\tau}_{n}(\{ y+\eta \nu_{y}\in R_{x}:\eta \in \tilde{A}
\}).\]

    By hypothesis $\nu_{n}$ converges vaguely to some nonzero measure $\nu$ with
support in $\{ x \}$,
    because $\nu_{n}(S^{\perp}_{x})=\tilde{\tau}_{n}(R_{x})\ge
    \tilde{\tau}_{n}(B(x,\delta))\ge\epsilon>0$. To get a contradiction it is enough to show that
    $\nu$ is dominated by a doubling measure in $S_{x}^{\bot}.$

    We define
    \[\gamma_{n}(A)=\tilde{\mu}_{n}(\{ y+\eta \nu_{y}\in R_{x}:\eta \in
\tilde{A}  \}),\quad
    A\subset S_{x}^{\perp},\]
    where as before $\tilde{\mu}_{n}$ is the ``projection'' of $\mu_{n}$ to
    $\R^{d}$.
    Observe that $\nu_{n}(A)\le \gamma_{n}(A)$ and that $\gamma_{n}$ are
doubling measures
    all with the same doubling constant. Indeed, for any $\delta>\alpha>0$
    there exist $y_{1},\dots ,y_{N}\in B(x,\delta)\cap S$ such that
    \[\{ y+\eta \nu_{y}\in R_{x}:|\eta|<2\alpha \}\subset
\bigcup_{j=1}^{N}B(y_{j},5\alpha/2),\quad
    \sum_{j=1}^{N}\chi_{B(y_{j},5\alpha/4)}\le C.\]
    The ``projection'' of the $\mu_{n}$ to $\R^{d}$
    are doubling measures all with the same doubling constant so
\begin{align*}
    \gamma_{n}(x+(-2\alpha,2\alpha)\nu_{x}) & \le
\sum_{j=1}^{N}\tilde{\mu}_{n}(B(y_{j},5\alpha/2))
    \\
    \le
    &
    C \sum_{j=1}^{N}\tilde{\mu}_{n}(B(y_{j},5\alpha/4))\le  C
\gamma_{n}(x+(\alpha,\alpha)\nu_{x}).
\end{align*}

    Therefore by (\ref{doubling}) we have $C,\gamma>0$ constants such that
    $\nu_{n}(x+(-r,r)\nu_{x})\le C r^{\gamma}$ and
    the same holds for $\nu$. Observe that $\supp \nu$
    has to be of Hausdorff dimension $\ge \gamma>0$ and this would contradict
    $\supp \nu=\{ x \}$.
\end{proof}

\section{Uniform norm case}\label{Secinfty}

    In this section we want to find sufficient conditions in the sequence
$\mathcal{E}=\{ E_{L} \}_{L\ge 0}$ in order to get
    the $L^{\infty}$-Logvinenko-Sereda property, i.e.
\begin{equation}
                   \label{infty}
    \sup_{u\in \S^{d}}|Q_{L}(u)|\le C \sup_{u\in E_{L}}|Q_{L}(u)|,\quad\mbox{for
any }Q_{L}\in \Pi_{L},
\end{equation}
    with $C$ a constant that does not depend on $L$.

    Our main result is the following:

\begin{theorem}
\label{lsinfty}
    If $\mathcal{E}$ is relatively dense, then it is
$L^{\infty}$-Logvinenko-Sereda.
\end{theorem}

\begin{remark} The converse is false because there exist discrete sets (so with
zero Lebesgue measure) with
    comparison property \eqref{infty}.
\end{remark}

    In \cite{MT} the authors deal with the weighted one-dimensional case of
Theorem~\ref{lsinfty}.
    In this uniform case it is a natural assumption to consider only bounded
weights.
    They considered the family of weights
    $\omega\ge 0$ such that
\begin{equation}\label{rh}
    \omega(u)\le \frac{C}{\sigma(B)}\int_{B}\omega(v) d\sigma(v),
\end{equation}
    for any spherical cap $B\subset \S^{d}$ and $u\in B.$
    Following \cite{CN} we call $RH_{\infty}$ this family.

\begin{definition}
    Let $\omega\ge 0$ be a function such that
    property \eqref{rh} holds for almost every $u\in \S^{d}$, we say that
$\omega$ is in
    the \emph{reverse H\"older class} $RH_{\infty}$.
\end{definition}

    To justify the name of this class observe that
    for $\omega\in RH_{\infty}$
    the reverse H\"older inequality
    \[\left( \frac{1}{\sigma(B)}\int_{B}\omega^{s}(u) d\sigma(u)\right)^{s}\le
\frac{C}{\sigma(B)}\int_{B}\omega(u) d\sigma(u),\quad B\subset
\S^{d} \mbox{ spherical cap}\]
    holds for each $s>1$, (i.e. $\omega \in RH_{s}$) and the best constant $C$
is bounded by the
    constant appearing in \eqref{rh}.
    And conversely, if the reverse H\"older inequality holds for each $s>1$
with a constant independent of $s$, then
    $\omega\in RH_{\infty}$, see \cite{CN}.

    Observe that $RH_{\infty} \subset A_{\infty}$.
    Roughly speaking $\omega$ belongs to $A_{1}$ if and only if $1/\omega \in
RH_{\infty}$. These weights can
    have high order zeros in $\S^{d}$.

    In this section we will prove the one-dimensional unweighted result first
and then
    extend it to $\S^{d}$.
    Using this unweighted case and adapting some results from \cite{Dai,MT} we
will prove the weighted
    result.

\begin{proof}
    We start with the one-dimensional case.
    Using the Lemma \ref{lemma} we get $h_{x}(E_{n})\ge \alpha$, for any $|x|=1-1/L$.
    Let $p$ be a polynomial of degree $L$, there exists an holomorphic
polynomial $q$ of degree
    $2L$ such that $|p|=|q|$ in $\S^{1}$.
    So for any $x\in \R^{2}$ with $|x|=1-1/L,$
\begin{align*}
    \log |q(x)| & \le h_{x}(E_{L})\log (\max_{E_{L}}|q|)+
    h_{x}(\S^{1}\setminus E_{L})\log (\max_{\S^{1}} |q|)
    \\
    & =
    \log \| q \|_{\S^{1}}+h_{x}(E_{L})\log \frac{\| q \|_{E_{L}}}{\| q
\|_{\S^{1}}}
    \le
    \log \| q \|_{\S^{1}}+\alpha\log \frac{\| q \|_{E_{L}}}{\| q \|_{\S^{1}}},
\end{align*}
    because $\| q \|_{E_{L}}/\| q \|_{\S^{1}}\le 1$ and so $|q(x)|\le \| p
\|_{E_{L}}^{\alpha}
    \| p \|_{\S^{1}}^{1-\alpha}$.
    Finally on can see that
    \[\max_{x\in \S^{1}}|q(x)|\le C \max_{|x|=1-1/L}|q(x)|,\]
    where $C$ is independent of $L$, see \cite[Lemma~2]{OS}.

    Now we consider the case $d>1.$
    Let $Q\in \Pi_{L}$ and suppose that
    $\max_{\S^{d}}|Q|=|Q(N)|=1$.
    We have that
\[\frac{\sigma( E_{L}\cap B(N,r/L))}{\sigma(B(N,r/L))}\ge \epsilon>0.\]
    Denoting $\tilde{\omega}=(\omega,0)\in \R^{d+1}$ for $\omega\in \S^{d-1}$ we
have
    that $G_{\omega}(\theta)=N \cos \theta+\tilde{\omega}\sin \theta$ is a
geodesic in $\S^{d}$
    if $\theta \in [-\pi,\pi]$. Therefore denoting $\S^{d-1}_{+}=\{ \omega\in
\S^{d-1}:\omega_{d}>0 \}$
\[\sigma(E_{L}\cap B(N,r/L))=\int_{\S^{d-1}_{+}}\int_{-r/L}^{r/L}
\chi_{E_{L}}(G_{\omega}(\theta))\sin^{d-1} \theta d\theta
d\omega.\]
    Now
\begin{eqnarray*}
\epsilon \sigma(B(N,r/L)) & \le & \int_{\S^{d-1}_{+}}
\int_{-r/L}^{r/L} \chi_{E_{L}}(G_{\omega}(\theta))\sin^{d-1}
\theta d\theta d\omega
\\
 & \le &
\int_{\S^{d-1}_{+}}\int_{-r/L}^{r/L}
\chi_{E_{L}}(G_{\omega}(\theta))\left( \frac{r}{L}\right)^{d-1}
d\theta d\omega
\\
& \le &
   \left( \frac{r}{L}\right)^{d-1} \int_{\S^{d-1}_{+}}  \sigma(E_{L}\cap
B_{\omega}(N,r/L)) d\omega
\end{eqnarray*}
    where $B_{\omega}(N,r/L)=\{ G_{\omega}(\theta):|\theta|\le r/L \}$.
    We get
$$\int_{\S^{d-1}_{+}}  \frac{\sigma (E_{L}\cap
B_{\omega}(N,r/L))}{\sigma(B_{\omega}(N,r/L))} d\omega \ge  C_{d}
\epsilon,$$
    and therefore there exists a direction $\omega\in \S^{d-1}_{+}$ such that
\begin{equation}
              \label{sep}
\frac{\sigma(E_{L}\cap B_{\omega}(N,r/L))}{\sigma
(B_{\omega}(N,r/L))}\ge C_{d}\epsilon >0.
\end{equation}
    We get the result applying the one-dimensional case to
    the trigonometric polynomial $Q(G_{\omega}(\theta) ).$
\end{proof}

    Using Theorem~\ref{lsinfty} we prove the following weighted
    version.

\begin{corollary}
 \label{lswinfty}
    If $\mathcal{E}$ is relatively dense and $\omega\in RH_{\infty}$, then
\begin{equation*}
                    \label{winfty}
    \sup_{u\in \S^{d}}|Q_{L}(u)|\omega(u)\le C \sup_{u\in
E_{L}}|Q_{L}(u)|\omega(u), \mbox{ for any }Q_{L}\in \Pi_{L},
\end{equation*}
    with $C$ a constant that does not depend on $L$.
\end{corollary}

\begin{remark} This result is optimal in some sense because there are unbounded
weights belonging to all reverse H\"older classes
    i.e. in particular $RH_{\infty}\varsubsetneq  \cap_{s>1} RH_{s}$, see
\cite[p.~2948]{CN}.
\end{remark}
\begin{proof}
    By definition of $RH_{\infty}$ weight
    \[\omega(u)\le C \omega_{L}(u)=\frac{1}{\sigma
(B(u,1/L))}\int_{B(u,1/L)}\omega(v) d\sigma(v).\]
    Now \cite[Lemma 4.6.]{Dai} provide us with $R_{L}\in \Pi_{L}$ nonnegative
such that
    for any $u\in \S^{d}$
    $\omega_{L}(u)\sim R_{L}(u)$,
    with constant depending only on the doubling constant for $\omega_{L}$.

    Now we want to construct a relatively dense regularization of
$E_{L}$
    that we will denote $E_{L}^{*}$.
    Given $\epsilon>0$ let $V=V_{\epsilon,L}\subset \S^{d}$ discrete and such
that
    \[\S^{d}\subset \bigcup_{v\in V}B(v,\epsilon/L), \mbox{ and }\sum_{v\in
V}\chi_{B(v,\epsilon/L)}(u)\le C_{d},\quad
    u\in \S^{d}.\]
    For $\delta>0,$ that we will determine afterwards,
    let \[V_{g}=\{ v\in V:\sigma(B(v,\epsilon/L)\cap E_{L})\ge \delta
\sigma(B(v,\epsilon/L))
    \}, \mbox{ and }E_{L}^{*}=\bigcup_{v\in V_{g}}B(v,\epsilon/L).\]
    We denote $V_{b}=V\setminus V_{g}$.
    Let $V(u)$ be the set of those $v\in V$ such that $B(v,\epsilon/L)\cap
B(u,r/L)\neq
    \emptyset$ and likewise we split $V(u)=V_{g}(u)\cup V_{b}(u)$
\begin{align}                                           \label{eq2}
    \sigma(B(u,r/2L)) & \le \sigma(\bigcup_{v\in
V_{g}(u)}B(v,\epsilon/L))+\sigma(\bigcup_{v\in
V_{b}(u)}B(v,\epsilon/L))
    \\
    &
    \le \sigma(E_{L}^{*}\cap B(u,r/L))+\sigma(\bigcup_{v\in
    V_{b}(u)}B(v,\epsilon/L)).\nonumber
\end{align}
    Using the relative density of $E_{L}$ and the property of being in $V_{b}$
we get
\begin{align*}
    \varrho \sigma(B(u,r/2L)) & \le \sigma(E_{L}\cap B(u,r/L))\\
    &\le
    C_{d}\delta \sigma(B(u,r/L))+\sigma\Bigl( E_{L}\cap \bigl( B(u,r/L)\setminus  \bigcup_{v\in
V_{b}(u)}B(v,\epsilon/L)\bigr)\Bigr),
\end{align*}
    so for $\delta$ small enough
\begin{equation}
\label{eq1}
    \sigma(B(u,r/L))- \sigma(\bigcup_{v\in
    V_{b}(u)}B(v,\epsilon/L))
    \ge
    \frac{\varrho}{2} \sigma(B(u,r/L)),
\end{equation}
    so using (\ref{eq1}) and (\ref{eq2}) we get
    \[ \frac{\varrho}{2}\sigma(B(u,r/2L))\le \sigma(E_{L}^{*}\cap B(u,r/L)),\]
    and thus $E_{L}^{*}$ is relatively dense.

Applying our unweighted result Theorem~\ref{lsinfty}
    to $E_{L}^{*}$ and to $M_{2L}\in \Pi_{2L},$ such that
    $M_{2L}=Q_{L}R_{L}$ in $\S^{d},$
    we get
    \[\sup_{u\in \S^{d}}|Q_{L}(u)|\omega(u)\lesssim \sup_{u\in
E_{L}^{*}}|Q_{L}(u)|\omega_{L}(u), \quad Q_{L}\in \Pi_{L}.\]

    We can take $\epsilon>0$ small enough so
    that spherical harmonics of degree $\le L$ are pointwise equivalents in
spherical caps of radius $\epsilon/L$ where they reach their
maximum.
    Indeed all we have to do is to apply Bernstein's inequality as we did in
proving Proposition~\ref{prop1}.

    Let $w \in B(v,\epsilon/L)$ with $v$ the center of a cap in $E_{L}^{*}.$ We
    apply the $A_{\infty}$ condition getting
\begin{align*}
    \omega_{L} & (w)=\frac{1}{\sigma
(B(w,1/L))}\int_{B(w,1/L)}\omega(u) d\sigma(u)
    \\
    &
    \le
    \frac{K}{\sigma (B(w,1/L))}\left(
\frac{\sigma(B(w,1/L))}{\sigma(B(v,\epsilon/L)\cap E_{L})}
    \right)^{s}\int_{B(v,\epsilon/L)\cap E_{L}}\omega(u) d\sigma(u)
    \\
    &
    \le
    \frac{C}{\sigma (B(w,1/L))}\left( \frac{\sigma(B(w,1/L))}{\delta
\sigma(B(v,\epsilon/L))}
    \right)^{s}\int_{B(v,\epsilon/L)\cap E_{L}}\omega(u) d\sigma(u)
    \\
    &
    =
    C_{\epsilon,\delta}L^{d}\int_{B(v,\epsilon/L)\cap E_{L}}\omega(u) d\sigma(u).
\end{align*}
    Finally, there exists $u\in V_{g}$ such that
    \[ \sup_{u\in
E_{L}^{*}}|Q_{L}(u)|\omega_{L}(u)=\sup_{u\in
B(v,\epsilon/L)}|Q_{L}(u)|\omega_{L}(u),  \]
    for any $w \in B(v,\epsilon/L)$
\begin{align*}
    \inf_{u\in B(v,\epsilon/L)}|Q_{L}(u)|\omega_{L}(w) & \le L^{d}\inf_{u\in
B(v,\epsilon/L)}|Q_{L}(u)|
    \int_{B(v,\epsilon/L)\cap E_{L}}\omega(z) d\sigma(z)
    \\
    &
    \le \sup_{u\in B(v,\epsilon/L)\cap E_{L}}|Q_{L}(u)|\omega(u),
\end{align*}
    and the result follows easily.
\end{proof}


\end{document}